\documentclass[a4paper, 12pt]{article}

\usepackage{amssymb,latexsym} \usepackage{mathtools}
\usepackage{graphicx,color,xcolor} \usepackage{hyperref}

\author{Johannes Sj\"ostrand\\
\small Centre de Math\'ematiques Laurent Schwartz\\ \small Ecole
Polytechnique\\ \small FR-91128 Palaiseau C\'edex\\ \footnotesize
johannes@math.polytechnique\\ \footnotesize and UMR 7640, CNRS }
\date{} \title{Some results on non-self-adjoint operators, a survey}

\newtheorem{dref}{Definition}[section] 
\newtheorem{theo}[dref]{Theorem} \newtheorem{theorem}[dref]{Theorem}
\newtheorem{prop}[dref]{Proposition} 
\newtheorem{ex}[dref]{Example}

\newcommand{\ekv}[2]{\begin{equation}\label{#1}#2\end{equation}}
\newcommand{\eekv}[3]{\begin{eqnarray}\label{#1}#2 \\ #3
\nonumber\end{eqnarray}}

\begin{document}
\maketitle
\begin{abstract} This text is a survey of recent results obtained by
the author and collaborators on different problems for
non-self-adjoint operators. The topics are: Kramers-Fokker-Planck type
operators, spectral asymptotics in two dimensions and Weyl asymptotics
for the eigenvalues of non-self-adjoint operators with small random
perturbations. In the introduction we also review the notion of
pseudo-spectrum and its relation to non-self-adjoint spectral
problems.
\end{abstract}
\tableofcontents
\section{Introduction}\label{int} \setcounter{equation}{0}
\subsection{Some background}\label{ld}
\par For self-adjoint and more generally normal operators on some
complex Hilbert space ${\cal H}$ we have a nice theory, including the
spectral theorem and a nice estimate on the norm of the resolvent:
\ekv{int.1} {\Vert (z-P)^{-1}\Vert \le ({\rm dist\,}(z,\sigma
(P)))^{-1},\quad \sigma (P)=\hbox{\, the spectrum of }P.}  This has a
consequence for the corresponding evolution problem: If $\sigma
(P)\subset\{ z\in{\bf C};\Re z\ge \lambda _0\}$, then \ekv{int.2}
{\Vert e^{-tP}\Vert \le e^{-\lambda _0t},\ t\ge 0 .}  However,
non-normal operators appear frequently: Scattering poles,
Convection-diffusion problems, Kramers-Fokker-Planck equation, damped
wave equations, linearized operators in fluid dynamics. Then
typically, $\Vert (z-P)^{-1}\Vert $ may be very large even when $z$ is
far from the spectrum. This implies mathematical difficulties:
\par\noindent -- When studying the distribution of eigenvalues,
\par\noindent -- When studying functions of the operator, like $
e^{-tP}$ and its norm.

\par It also implies numerical difficulties like:
\par\noindent -- Eigenvalue instability.

\par There are (in the author's opinion) two ways out:
\par\noindent -- Change the Hilbert space norm to make the operators
look more normal. (Complex scaling methods.)
\par\noindent -- Recognize that the region of the $z$-plane where
$\Vert (z-P)^{-1}\Vert $ is large, has its own
interest. (Pseudospectrum.)

\par The option to choose depends on the problem. \begin{itemize}
\item In some problems, like those related to scattering poles, there
is no obvious choice of Hilbert space and we are free to make the most
natural one. This option is particularly natural when considering a
differential operator with analytic coefficients.
\item In other problems the canonical Hilbert space is $L^2$ and we
are at most allowed to change the norm into an equivalent one. Here
the notion of pseudospectrum is likely to be important.
\end{itemize}

 Let $P:{\cal H}\to{\cal H}$ be closed, densely defined, ${\cal H}$ a
complex Hilbert space and let $\rho (P)={\bf C}\setminus \sigma (P)$
denote the resolvent set. The notion of pseudospectrum is important in
numerical analysis and we refer to L.N.~Trefethen \cite{Tr},
Trefethen--M.~Embree \cite{TrEm} and further references given
there. Thanks to works of E.B.~Davies \cite{Da}, \cite{Da2},
M.~Zworski \cite{Zw} and others it has become popular in the
non-self-adjoint spectral theory of differential operators.
\begin{dref} \label{int1} Let $\epsilon >0$. The $\epsilon
$-pseudo\-spectrum of $P$ is
$$
\sigma _\epsilon (P):=\sigma (P)\cup\{ z\in \rho (P);\, \Vert
(z-P)^{-1}\Vert > 1/\epsilon \}.
$$ \end{dref}
Unlike the spectrum, the pseudospectrum will in general change when we
change the norm on ${\cal H}$. Moreover, it can be characterized as a
set of spectral instability as follows from the following version of a
theorem of Roch-Silberman \cite{RoSi}:
\begin{theo}\label{int2}\end{theo}
$$
\sigma _\epsilon (P)=\bigcup_{Q\in{\cal L}({\cal H},{\cal H})\atop
\Vert Q\Vert < \epsilon }\sigma (P+Q).
$$

In his survey \cite{Tr} L.N.~Trefethen discusses some linearized
operators from fluid dynamics:
\begin{itemize}
\item Orr-Sommerfeld equation ({\small Orzag, Reddy, Schmid,
Hennigson}).
\item Plane Poiseuille flow ({\small L.N and A.N Trefethen, Schmid}).
\item Pipe Poiseuille flow ({\small L.N and A.N Trefethen, Reddy,
Driscoll}),
\end{itemize} and to what extent stability can be predicted from the
sudy of the spectrum of these non-self-adjoint operators: Eigenvalue
analysis alone leads in some cases to the prediction of stability for
Reynolds numbers $R<5772$. Experimentally however, we have stability
only for $R<1000$.

\par The rough explanation of this is that the $\epsilon
$-pseudospectrum (for a suitable $\epsilon $) crosses the imaginary
axis before the spectrum does, when $R$ increases. Then $\Vert
e^{-tP}\Vert $ will grow fast for a limited time even though the
growth for very large times is determined by the spectrum. However,
since $P$ appears as a linearization of a non-linear problem, that
suffices to cause instability.

\par In the case of differential operators the pseudospectral
phenomenon is very general and related to classical works in PDE on
local solvability and non-hypoellipticity.  E.B.~Davies \cite{Da}
studied the non-self-adjoint semiclasscial Schr\"odinger operator with
a smooth (complex-valued) potential in dimension 1 and showed under
``generic'' assumptions that one can construct quasimodes with the
spectral parameter varying in an open complex set, containing points
that are possibly very far from the spectrum (as can be verified in
the case of the complex harmonic oscillator).  M.~Zworski \cite{Zw}
observed that this is essentially a rediscovery of an old result of
H\"ormander \cite{Ho2, Ho3}, and was able to generalize considerably
Davies' result by adapting the one of H\"ormander to the
semi-classical case. With N.~Dencker and M.Zworski \cite{DeSjZw} we
also gave a direct proof and a corresponding adaptation of old results
of Sato-Kawai-Kashiwara \cite{SKK} to the analytic case:
\begin{theo}\label{int3}(\cite{Zw}, \cite{DeSjZw}) Let \ekv{int.3} {
P(x,hD_x)=\sum_{\vert \alpha \vert \le m}a_\alpha (x)(hD_x)^\alpha ,\
D_x={\partial \over\partial x} } have smooth coefficients in the open
set $\Omega \subset{\bf R}^n$. Put $p(x,\xi )=\sum_{\vert \alpha \vert
\le m}a_\alpha (x)\xi ^\alpha $. Assume $z=p(x_0,\xi _0)$ with the
Poisson bracket ${1\over i}\{ p,\overline{p}\}(x_0,\xi_0) >0$. Then
$\exists$ $u=u_h\in C_0^\infty (\Omega )$, with $\Vert u\Vert=1$,
$\Vert (P-z)u\Vert ={\cal O}(h^\infty )$, when $h\to 0$.
\par\noindent Analytic case: Can replace "$h^\infty $" by
"$e^{-1/Ch}$".\end{theo}

Here, we have used standard multi-index notation: $\xi ^\alpha =\xi
_1^{\alpha _1}\cdot ...\cdot \xi _n^{\alpha _n}$, $|\alpha |=\alpha
_1+...+\alpha _n$, and the norms $\Vert \cdot \Vert$ are the ones of
$L^2$ or $\ell ^2$ if nothing else is indicated.

This result was subsequently generalized by K.\ Pravda-Starov,
\cite{Pr}.  Notice that this implies that when the theorem applies and
if the resolvent $(P-z)^{-1}$ exists then its norm is greater than any
negative power of $h$ when $h\to 0$.
\begin{ex}\label{in4} Let $P=-h^2\Delta +V(x)$, $p(x,\xi
)=\xi^2+V(x)$, $\xi ^2=\xi _1^2+...+\xi _n^2$. Then ${1\over
i}\{p,\overline{p}\}=-4\xi \cdot \Im V'(x)$.
\end{ex}

\par ``Generically'', if $z=p(x,\xi )$, then $\{ p,\overline{p}\}
(x,\xi )\ne 0$ and one can show quite generally that if this happens
then there is also another point $(y,\eta )$ with $p(y,\eta )=z$ such
that $\{ p,\overline{p}\} (x,\xi )$ takes the opposite sign.  This
justifies the following simplified terminology in the semi-clas\-sical
limit:

\par The semi-classical pseudospectrum of $P$ is the range ${\cal
R}(p)$ of $p$.

\par In \cite{DeSjZw} we also showed under suitable assumptions
(inspired from scattering theory and from the theory of sub-elliptic
operators), that there is no spectrum near the boundary of the
semi-classical pseudospectrum and that we may have quite a good
control of the norm of the resolvent there. Generalizations to the
case of systems were given by Dencker \cite{De}.

\subsection{The topics of this survey}

We will discuss three subjects involving non-self-adjoint differential
and pseudodifferential operators, We will always wok in the
semi-classical limit, which means that our operators are of the form
$P(x,hD_x;h)$, where $P$ is a suitable symbol and $0<h\ll 1$. It is
quite clear however that some of our results will also apply to
non-semiclassical situations in the limit of large eigenvalues.

\par The subjects are:
\begin{itemize}
\item The Kramers-Fokker-Planck operator,
\item Bohr-Sommerfeld rules in dimension 2,
\item Weyl asymptotics for non-self-adjoint operators with small
random perturbations.
\end{itemize} and most of the works discussed are the results of
collaborations with A.~Melin, M.~Hitrik, F.~H\'erau, C.Stolk, S.~V{\~
u} Ng{\d o}c and M.~Hager.

\par In the first two topics we exploit the possibility of changing
the Hilbert space norm by introducing exponential weights on phase
space.

\par In the case of the Kramers-Fokker-Planck operator, we make no
analyticity assumptions and the phase-space weights are
correspondingly quite weak. In this case however our operator is a
differential one, so we are allowed to apply strong exponential
weights depending only on the base variables, and this is important
when studying small exponential corrections of the eigenvalues via the
so called tunnel effect.

\par For the Bohr-Sommerfeld rules, we make analyticity assumptions
that allow stronger phase space weights. In both cases the effect of
the exponential weights is to make the operator under consideration
more normal.

\par In the third topic, we do not use any deformations of the given
Hilbert space, but exponential weights play an important role at
another level, namely to count zeros of holomorphic functions with
exponential growth.

\par The pseudospectrum will not be discussed explicitly below. In the
Kramers-Fokker-Planck case, the problems are located near the boundary
of the semi-classical pseudospectrum, and it turns out that we have a
very nice control of the resolvent there. In the 2 dimensional
Bohr-Sommerfeld rules, we have stronger exponential weights,
reflecting stronger pseudospectral phenomena. Finally in the subject
of Weyl asymptotics, we often have strong pseudospectral behaviour for
the unperturbed operator.  From the proofs it appears that the random
perturbations will weaken the pseudospectral behaviour and this might
have very interesting consequences for the associated evolution
problems. This is still very much an open problem.

\section{Kramers-Fokker-Planck type operators, spectrum and return to
equilibrium}\label{kfp} \setcounter{equation}{0}
\subsection{Introduction}
\par There has been a renewed interest in the problem of {``return to
equilibrium''} for various 2nd order operators.  One example is the
Kramers-Fokker-Planck operator: \ekv{0.1} {P= y\cdot h\partial
_x-V'(x)\cdot h\partial _y+{\gamma \over 2}(-h\partial _y+y) \cdot
(h\partial _y+y), } where $x,y\in {\bf R}^n$ correspond respectively
to position and speed of the particles and $h>0$ corresponds to
temperature. The constant $\gamma >0$ is the friction. (Since we will
only discuss $L^2$ aspects we here present right away an adapted
version of the operator, obtained after conjugation by a Maxwellian
factor.)

\par The associated evolution equation is:
$$
(h\partial _t+P)u(t,x,y)=0.
$$
{\it Problem of return to equilibrium:} Study the rate of convergence
of $u(t,x,y)$ to a multiple of the ``ground state''
$u_0(x,y)=e^{-(y^2/2+V(x))/h}$ when $t\to +\infty $, assuming that
$V(x)\to +\infty $ sufficiently fast when $x\to \infty $ so that
$u_0\in L^2({\bf R}^{2n})$. Notice here that $P(u_0)=0$ and that the
vector field part of $P$ is $h$ times the Hamilton field of
$y^2/2+V(x)$, when we identify ${\bf R}_{x,y}^{2n}$ with the cotangent
space of ${\bf R}^n_x$.

\par A closely related problem is to study the difference between the
first eigenvalue (0) and the next one, $\mu (h)$. (Since our operator
is non-self-adjoint, this is only a very approximate formulation
however.)

Some contributions: L.~Desvillettes--C.~Villani \cite{DeVi},
J.P.~Eckmann--M.~Hairer \cite{EcHa}, F.~H\'erau--F.~Nier \cite{HeNi},
B.~Helffer--F.~Nier \cite{HelNi}, Villani \cite{Vi}. In the work
\cite{HeNi} precise estimates on the exponential rates of return to
equilibrium were obtained with methods close to those used in
hypoellipticity studies and this work was our starting point. With
H\'erau and C.Stolk \cite{HeSjSt} we made a study in the
semi-classical limit and studied small eigenvalues modulo ${\cal
O}(h^{\infty })$. More recently with H\'erau and M.~Hitrik
\cite{HeHiSj} we have made a precise study of the exponential decay of
$\mu (h)$ when $V$ has two local minima (and in that case $\mu (h)$
turns out to be real). This involves tunneling, i.e. the study of the
exponential decay of eigenfunctions. As an application we have a
precise result on the return to equilibrium \cite{HeHiSj2}.  This has
many similarities with older work on the tunnel effect for
Schr\"odinger operators in the semi-classical limit by
B.~Helffer--Sj\"ostrand \cite{HelSj1, HelSj2} and B.~Simon \cite{Si}
but for the Kramers-Fokker-Planck operator the problem is richer and
more difficult since $P$ is neither elliptic nor self-adjoint. We have
used a supersymmetry observation of J.M.~Bismut \cite{Bi} and
J.~Tailleur--S.~Tanase-Nicola--J.~Kurchan \cite{TaTaKu}, allowing
arguments similar to those for the standard Witten complex
\cite{HelSj2}.
\subsection{Statement of the main results}\label{kfpst}

 Let $P$ be given by (\ref{0.1}) where $V\in C^\infty ({\bf R}^n;{\bf
R})$, and \ekv{1.1} {\partial ^\alpha V(x)={\cal O}(1),\ \vert \alpha
\vert\ge 2,} \ekv{1.2} { \vert \nabla V(x) \vert\ge 1/C,\ \vert x
\vert \ge C, } \ekv{1.3} { V\hbox{ is a Morse function.}  } We also
let $P$ denote the graph closure of $P$ from ${\cal S}({\bf R}^{2n})$
which coincides with the maximal extension of $P$ in $L^2$ (see
\cite{HeNi, HelNi, HeHiSj2}).  We have $\Re P\ge 0$ and the spectrum
of $P$ is contained in the right half plane. In \cite{HeSjSt} the
spectrum in any strip $0\le \Re z\le Ch$ (and actually in a larger
parabolic neighborhood of the imaginary axis, in the spirit of
\cite{HeNi}) was determined asymptotically ${\rm mod\,}({\cal
O}(h^\infty ))$. It is discrete and contained in a sector $\vert \Im z
\vert\le C\Re z+{\cal O}(h^\infty )$:
\begin{theorem}\label{kfp1} The eigenvalues in the strip $0\le \Re
z\le Ch$ are of the form \ekv{1.4} { \lambda
_{j,k}(h)\sim h(\mu _{j,k}+h^{1/N_{j,k}}\mu _{j,k,1}+ h^{2/N_{j,k}}\mu
_{j,k,2}+..)  } where $\mu _{j,k}$ are the eigenvalues of the
quadratic approximation (``non-selfadjoint oscillator'')
$$
y\cdot \partial _x-V''(x_j)x\cdot \partial _y+{\gamma \over
2}(-\partial _y+y)\cdot (\partial _y+y),
$$
at the points $(x_j,0)$, where $x_j$ are the critical points of $V$.
\end{theorem}

 The $\mu _{j,k}$ are known explicitly and it follows that when $x_j$
is not a local minimum, then $\Re\lambda _{j,k}\ge h/C$ for some
$C>0$.  When $x_j$ is a local minimum, then precisely one of the
$\lambda_{j,k}$ is ${\cal O}(h^\infty )$ while the others have real
part $\ge h/C$. Furthermore, when $V\to +\infty $ as $x\to \infty $,
then $0$ is a simple eigenvalue. {\it In particular, if $V$ has only
one local minimum, then
$$
\inf \Re (\sigma (P)\setminus \{ 0\})\sim h(\mu _1+h\mu _2+\dots
),\quad \mu _1>0.
$$
}(or possibly an expansion in fractional powers) and we obtained a
corresponding result for the problem of return to equilibrium. It
should be added that when $\mu _{j,k}$ is a simple eigenvalue of the
quadratic approximation then $N_{j,k}=1$ so there are no fractional
powers of $h$ in (\ref{1.4}).

\par The following is the main new result that we obtained with
F.~H\'erau and M.~Hitrik in \cite{HeHiSj}:
\begin{theorem}\label{kfp2} Assume that $V$ has precisely 3 critical
points; 2 local minima, $x_{\pm 1}$ and one ``saddle point'', $x_0$ of
index 1. Then for $C>0$ sufficiently large and $h$ sufficiently small,
$P$ has precisely 2 eigenvalues in the strip $0\le \Re z\le h/C$,
namely $0$ and $\mu (h)$, where $\mu (h)$ is real and of the form
\ekv{1.5} { \mu (h)=h(a_1(h)e^{-2S_1/h}+a_{-1}(h)e^{-2S_{-1}/h}), }
where $a_j$ are real,
$$
a_j(h)\sim a_{j,0}+ha_{j,1}+...,\ h\to 0,\quad a_{j,0}>0,
$$
$$S_j=V(x_0)-V(x_j).$$\end{theorem}

\par As for the problem of return to equilibrium, we obtained the
following result with F.~H\'erau and M.~Hitrik in \cite{HeHiSj2}:
\begin{theo}\label{kfp2.5} We make the same assumptions as in Theorem
\ref{kfp2} and let $\Pi_j$ be the spectral projection associated with
the eigenvalue $\mu_j$, $j=0,1$, where $\mu _0=0$, $\mu _1=\mu
(h)$. Then we have
\begin{equation}
\label{eq04} \Pi_j={\cal O}(1):L^2\to L^2,\quad h\rightarrow 0.
\end{equation} We have furthermore, uniformly as $t\geq 0$ and
$h\rightarrow 0$,
\begin{equation}
\label{eq05} e^{-tP/h}=\Pi_0+e^{-t\mu_1/h}\Pi_1+{\cal O}(1) e^{-t/C},
\mbox{ in }{\cal L}(L^2,L^2),
\end{equation} where $C>0$ is a constant.
\end{theo}

\par Actually, as we shall see in the outline of the proofs, these
results (as well as (\ref{1.4})) hold for more general classes of
supersymmetric operators.

\subsection{A partial generalization of \cite{HeSjSt}}\label{kfpgen}

Consider on ${\bf R}^n$ ($2n$ is now replaced by $n$):
\begin{eqnarray*} P&=&\sum_{j,k}hD_{x_j}b_{j,k}(x)hD_{x_k}+\\
&&{1\over 2}\sum_j(c_j(x)h\partial _{x_j} +h\partial _{x_j}\circ
c_j(x))+p_0(x)\\ &=&P_2+iP_1+P_0,
\end{eqnarray*} where $b_{j,k}, c_j, p_0$ are real and smooth.  The
associated symbols are:
\begin{eqnarray*} p(x,\xi )&=&p_2(x,\xi )+ip_1(x,\xi )+p_0(x),\\
p_2&=&\sum b_{j,k}\xi _j\xi _k,\ p_1=\sum c_j\xi _j.
\end{eqnarray*} Assume,
$$
p_2\ge 0,\ p_0\ge 0,
$$
\begin{eqnarray*}
\partial _x^\alpha b_{j,k}&=&{\cal O}(1),\ \vert \alpha \vert\ge 0,\\
\partial _x^\alpha c_j&=&{\cal O}(1),\ \vert \alpha \vert\ge 1,\\
\partial _x^\alpha p_0&=&{\cal O}(1),\ \vert \alpha \vert\ge 2.
\end{eqnarray*} {Assume that} $$\{x;\,p_0(x)=c_1(x)=..=c_n(x)=0\}$$ is
finite $=\{ x_1,...,x_N\}$ and put ${\cal C}=\{ \rho _1,...,\rho
_n\}$, $\rho _j=(x_j,0)$.  {Put}
$$\widetilde{p}(x,\xi )=\langle \xi  \rangle^{-2}
p_2(x,\xi )+p_0(x),\ \langle \xi \rangle =\sqrt{1+\vert \xi \vert
^2}$$
$$
\langle \widetilde{p} \rangle_{T_0}={1\over T_0}\int_{-T_0/2}^{T_0/2}
\widetilde{p} \circ \exp (tH_{p_1})dt,\ T_0>0 \hbox{ fixed.}
$$
Here in general we let $H_a=a'_\xi \cdot \frac{\partial }{\partial x}
-a'_x \cdot \frac{\partial }{\partial \xi } $ denote the Hamilton
field of the $C^1$-function $a=a(x,\xi )$.

\par Dynamical assumptions: Near each $\rho _j$ we have $\langle
\widetilde{p} \rangle_{T_0}\sim \vert \rho -\rho _j \vert^2$ and in
any compact set disjoint from ${\cal C}$ we have $\langle
\widetilde{p} \rangle_{T_0}\ge 1/C$. (Near infinity this last
assumption has to be modified slightly and we refer to \cite{HeHiSj}
for the details.) The following result from \cite{HeHiSj} is very
close to the main result of \cite{HeSjSt} and generalizes Theorem
\ref{kfp1}:

\begin{theorem} \label{kfp3} Under the above assumptions, the spectrum
of $P$ is discrete in any band $0\le \Re z\le Ch$ and the eigenvalues
have asymptotic expansions as in (\ref{1.4}).\end{theorem}

\par Put $$q(x,\xi )=-p(x,i\xi )=p_2(x,\xi )+p_1(x,\xi )-p_0(x).$$

\par The linearization of the Hamilton field $H_q$ at $\rho _j$ (for
any fixed $j$) has eigenvalues $\pm \alpha _k$, $k=1,..,n$ with real
part $\ne 0$.  Let $\Lambda _+=\Lambda _{+,j}$ be the unstable
manifold through $\rho _j$ for the $H_q$-flow. Then $\Lambda _+$ is
Lagrangian and of the form $\xi =\phi '_+(x)$ near $x_j$ ($\phi_+
=\phi _{+,j}$), where
$$\phi_+ (x_j)=0,\ \phi_+' (x_j)=0,\ 
\phi_+'' (x_j)>0.$$ The next result is from \cite{HeHiSj}:
\begin{theorem}\label{kfp4} Let $\lambda _{j,k}(h)$ be a simple
eigenvalue as in (\ref{1.4}) and assume there is no other eigenvalue
in a disc $D(\lambda _{j,k},h/C)$ for some $C>0$.  Then, in the $L^2$ sense, the
corresponding eigenfunction is of the form $e^{-\phi
_+(x)/h}(a(x;h)+{\cal O}(h^\infty ))$ near $x_j$, where $a(x;h)$ is
smooth in $x$ with an asymptotic expansion in powers of $h$. Away from
a small neighborhood of $x_j$ it is exponentially decreasing.
\end{theorem}

The proof of the first theorem uses microlocal weak exponential
estimates, while the one of the last theorem also uses local
exponential estimates.

\subsection{Averaging and exponential weights.}\label{kfpav}

The basic idea of the proof of Theorem \ref{kfp3} is taken from
\cite{HeSjSt}, but we reworked it in order to allow for
non-hypoelliptic operators. We will introduce a weight on $T^*{\bf
R}^n$ of the form \ekv{pr.1} { \psi _\epsilon =-\int
J(\frac{t}{T_0})\widetilde{p}_\epsilon \circ \exp (tH_{p_1}) dt, } for
$0<\epsilon \ll 1$. Here $J(t)$ is the {odd} function given by
\ekv{pr.2} { J(t)= \left\{ \begin{array}{ll} & 0,\ \vert t\vert \ge
\frac{1}{2}, \\ & \frac{1}{2}-t,\ 0<t\le \frac{1}{2},
\end{array} \right.  } and we choose $\widetilde{p}_\epsilon (\rho )$
to be equal to $\widetilde{p}(\rho )$ when ${\rm dist\,}(\rho ,{\cal
C})\le \epsilon $, and flatten out to $\epsilon \widetilde{p}$ away
from a fixed neighborhood of ${\cal C}$ in such a way that
$\widetilde{p}_\epsilon ={\cal O}(\epsilon )$.  Then \ekv{pr.3} {
H_{p_1}\psi _\epsilon =\langle \widetilde{p}_\epsilon \rangle _{T_0}
-\widetilde{p}_\epsilon .  }

We let $\epsilon =Ah$ where $A\gg 1$ is independent of $h$. {Then the
weight $\exp (\psi _\epsilon /h)$ is uniformly bounded when $h\to 0$.}
Indeed, $\psi _\epsilon ={\cal O}(h)$.

\par Using Fourier integral operators with complex phase, we can
define a Hilbert space of functions that are ``microlocally ${\cal
O}(\exp (\psi _\epsilon /h))$ in the $L^2$ sense''. The norm is
uniformly equivalent to the one of $L^2$, but the natural leading
symbol of $P$, acting in the new space, becomes { \ekv{pr.4}{p(\exp
(iH_{\psi_\epsilon })(\rho )),\ \rho \in T^*{\bf R}^n} }which by
Taylor expansion has {real part $\approx p_2(\rho )+p_0(\rho )+\langle
\widetilde{p}_\epsilon \rangle-\widetilde{p}_\epsilon $.}
\par Very roughly, {the real part of the new symbol is $\ge \epsilon $
away from ${\cal C}$ and behaves like ${\rm dist\,}(\rho ,{\cal C})^2$
in a $\sqrt{\epsilon }$-neighborhood of ${\cal C}$.} This can be used
to show that the spectrum of $P$ (viewed as an operator on the
weighted space) in a band $0\le \Re z <\epsilon /C $ comes from an
$\sqrt{\epsilon }$-neighborhood of ${\cal C}$. In such a neighborhood,
we can treat $P$ as an elliptic operator and the spectrum is to
leading order determined by the quadratic approximation of the dilated
symbol (\ref{pr.4}). This gives Theorem \ref{kfp3}.

We next turn to the proof of Theorem \ref{kfp4}, and we work near a
point $\rho _j=(x_j,\xi _j)\in {\cal C}$. Recall that $\Lambda _+:\xi
=\phi _+'(x)$ is the unstable manifold for the $H_q$-flow, where
$q(x,\xi )=-p(x,i\xi )$. We have {$q(x,\phi _+'(x))=0$}.

\par In general, if $\psi \in C^\infty $ is real, then $P_\psi
:=e^{\psi /h}\circ P\circ e^{-\psi /h}$ has the symbol \ekv{pr.5} {
p_\psi (x,\xi )=p_2(x,\xi ){-q(x,\psi '(x))}+i(q'_\xi (x,\psi
'(x))\cdot \xi }

\begin{itemize}
\item As long as $q(x,\psi '(x))\le 0$, we have $\Re p_\psi \ge 0$ and
we may hope to establish good apriori estimates for $P_\psi $.
\item This is the case for $\psi =0$ and for $\psi =\phi _+$. Using
the convexity of $q(x,\cdot )$, we get suitable weights $\psi $ with
$q(x,\psi '(x))\le 0$, equal to $\phi _+(x)$ near $x_j$, strictly
positive away from $x_j$ and constant outside a neighborhood of that
point.
\item It follows that the eigenfunction in Theorem \ref{kfp4} is
(roughly) ${\cal O}(e^{-\phi _+(x)/h})$ near $x_j$ in the $L^2$ sense.
\item On the other hand, we have quasi-modes of the form
$a(x;h)e^{-\phi _+(x)/h}$ as in \cite{HelSj1}.
\item Applying the exponentially weighted estimates, indicated above,
to the difference of the eigenfunction and the quasi-mode, we then get
Theorem \ref{kfp4}.
\end{itemize}

\subsection{Supersymmetry and the proof of Theorem
\ref{kfp2}}\label{kfpss}

\par We review the supersymmetry from \cite{Bi}, \cite{TaTaKu}, see
also G.~Lebeau \cite{Le}.  Let $A(x):T_x^*{\bf R}^n\to T_x{\bf R}^n$
be linear, invertible and smooth in $x$. Then we have the
nondegenerate bilinear form
$$
\langle u\vert v \rangle_{A(x)}=\langle \wedge^kA(x)u|v \rangle,\
u,v\in\wedge^kT_x^*{\bf R}^n,
$$
and we also write $(u|v)_{A(x)}=\langle u|\overline{v}
\rangle_{A(x)}$.

\par If $u,v$ are smooth $k$-forms with compact support, put
$$
(u|v)_A=\int (u(x)|v(x))_{A(x)}dx.
$$
The formal ``adjoint'' $Q^{A,*}$ of an operator $Q$ is then given by
$$
(Qu|v)_A=(u|Q^{A,*}v)_A.
$$

 Let $\phi :{\bf R}^n\to {\bf R}$ be a smooth Morse function with
$\partial ^\alpha \phi $ bounded for $\vert \alpha \vert\ge 2$ and
with $|\nabla \phi |\ge 1/C$ for $|x|\ge C$. Introduce the
Witten-De Rham complex:
$$
d_\phi =e^{-{\phi \over h}}\circ hd\circ e^{{\phi \over h}}= \sum_j
(h\partial _{x_j}+\partial _{x_j}\phi )\circ dx_j^\wedge ,
$$
where $d$ denotes exterior differentiation and $dx_j^\wedge$ left
exterior multiplication with $dx_j$.  The corresponding Laplacian is
then: $-\Delta_A=d_\phi ^{A,*}d_\phi +d_\phi d_\phi^{A,*}$.  Its
restriction to $q$-forms will be denoted by $-\Delta _A^{(q)}$. Notice
that:
$$-\Delta ^{(0)}_A(e^{-\phi /h})=0.$$

Write $A=B+C$ with $B^t=B$, $C^t=-C$. $-\Delta_A$ is a second order
differential operator with scalar principal symbol in the
semi-classical sense (${h\over i}{\partial \over \partial x_j}\mapsto
\xi _j$) of the form:
$$p(x,\xi )=\sum_{j,k}b_{j,k}(\xi _j\xi _k+\partial _{x_j}\phi
\partial _{x_k}\phi )+2i\sum_{j,k}c_{j,k}\partial _{x_k}\phi \,\xi _j.
$$
\par\noindent {\bf Example.} Replace $n$ by $2n$, $x$ by $(x,y)$, let
$$
A={1\over 2}\left(\begin{array}{ccc}0 &I\\ -I
&\gamma \end{array}\right).
$$
Then
\begin{eqnarray*} -\Delta _A^{(0)}&=&h(\phi '_y\cdot \partial _x-\phi
'_x\cdot \partial _y) \\ &&+ {\gamma \over 2}\sum_{j}(-h\partial
_{y_j}+\partial _{y_j}\phi ) (h\partial _{y_j}+\partial _{y_j}\phi ).
\end{eqnarray*} When $\phi ={y^2/2}+V(x)$ we recover the KFP operator
(\ref{0.1})

\par The results of Subsection \ref{kfpgen} apply, if we make the
additional dynamical assumptions there; $-\Delta ^{(q)}_A$ has an
asymptotic eigenvalue $=o(h)$ associated to the critical point $x_j$
precisely when the index of $x_j$ is equal to $q$ (as for the Witten
complex and analogous complexes in several complex variables). In
order to cover the cases $q>0$ we also assume that \ekv{3.1}{A={\rm
Const.}}  {\it The Double well case.}  Keep the assumption
(\ref{3.1}).  Assume that $\phi $ is a Morse function with $|\nabla
\phi |\ge 1/C$ for $|x|\ge C$ such that $-\Delta _A$ satisfies the
extra dynamical conditions of Subsection \ref{kfpgen} and having
precisely three critical points, two local minima $U_{\pm 1}$ and a
saddle point $U_0$ of index 1.

\par Then $-\Delta _A^{(0)}$ has precisely 2 eigenvalues: $0,\, \mu $
that are $o(h)$ while $-\Delta _A^{(1)}$ has precisely one such
eigenvalue: $\mu $. (Here we use as in the study of the Witten
complex, that $d_\phi $ and $d_{\phi }^{A,*}$ intertwine our
Laplacians in degeree 0 and 1. The detailed justification is more
complicated however.)  $e^{-\phi /h}$ is the eigenfunction of $\Delta
_A^{(0)}$ corresponding to the eigenvalue 0. Let $S_j=\phi (U_0)-\phi
(U_j)$, $j=\pm 1$, and let $D_j$ be the connected component of $\{
x\in {\bf R}^n;\, \phi (x) < \phi (U_0)\}$ containing $U_j$ in its
interior.

\par

\par Let $E^{(q)}$ be the corresponding spectral subspaces so that
${\rm dim\,}E^{(0)}=2$, ${\rm dim\,}E^{(1)}=1$. Truncated versions of
the function $e^{-\phi (x)/h}$ can be used as approximate
eigenfunctions, and we can show:
\begin{prop}\label{kfp5} $E^{(0)}$ has a basis $e_1,e_{-1}$, where
$$e_j=\chi_j(x)e^{-{1\over h}(\phi (x)-\phi (U_j))}+{\cal
O}( e^{-{1\over h}(S_j-\epsilon )}),\mbox{ in the $L^2$-sense.} $$ Here,
we let $\chi _j\in C_0^\infty (D_j)$ be equal to 1 on $\{ x\in D_j;\,
\phi (x)\le \phi (U_0)-\epsilon \}$.\end{prop} The theorems
\ref{kfp3}, \ref{kfp4} can be adapted to $-\Delta _A^{(1)}$ and lead
to:
\begin{prop}\label{kfp6} $E^{(1)}={\bf C}e_0$, where
$$e_0(x)=\chi _0(x)a_0(x;h)e^{-{1\over h}\phi _+(x)}+{\cal O}
(e^{-\epsilon _0/h}),$$ $\phi _+(x)\sim (x-U_0)^2,$ $\epsilon _0>0$ is
small enough, $a_0$ is an elliptic symbol, $\chi_0\in C_0^\infty ({\bf
R}^n)$, $\chi_0=1$ near $U_0$.
\end{prop}

\par Let the matrices of $d_\phi :E^{(0)}\to E^{(1)}$ and
$d_\phi^{A,*} :E^{(1)}\to E^{(0)}$ with respect to the bases $\{
e_{-1}, e_1\}$ and $\{ e_0\}$ be
$$
\left(\begin{array}{ccc} \lambda_{-1} & \lambda_{1}\end{array}\right)
\hbox{ and } \left (\begin{array}{ccc}\lambda _{-1}^* \\ \lambda
_1^*\end{array}\right)\hbox{ respectively}.
$$
Using the preceding two results in the spirit of tunneling estimates
and computations of Helffer--Sj\"ostrand (\cite{HelSj1, HelSj2}) we
can show:

\begin{prop}\label{kfp7} Put $S_j=\phi (U_0)-\phi (U_j)$, $j=\pm
1$. Then we have
$$
\left(\begin{array}{ccc} \lambda _{-1}\\ \lambda _1\end{array}\right)
=h^{1\over 2}(I+{\cal O}(e^{-{1\over Ch}})) \left(\begin{array}{ccc}
\ell_{-1}(h)e^{-S_{-1}/h}\\ \ell_1(h)e^{-S_1/h}
\end{array}\right),
$$

$$
\left(\begin{array}{ccc} \lambda ^*_{-1}\\ \lambda
_1^*\end{array}\right)=h^{1\over 2}(I+{\cal O}(e^{-{1\over Ch}}))
\left(\begin{array}{ccc} \ell_{-1}^*(h)e^{-S_{-1}/h}\\
\ell_1^*(h)e^{-S_1/h}
\end{array}\right),
$$
where $\ell_{\pm 1}$, $\ell_{\pm 1}^*$ are real elliptic symbols of
order $0$ such that $\ell_{j}\ell_j^*>0$, $j=\pm 1$.
\end{prop}

\par From this we get Theorem \ref{1.2}, since $\mu =\lambda
^*_{-1}\lambda _{-1}+ \lambda _1^*\lambda
_1$.\hfill{$\square$}\medskip

Thanks to the fact that we have only two local minima, certain
simplifications were possible in the proof. In particular it was
sufficent to control the exponential decay of general eigenfunctions
in some small neighborhood of the critical points. For more general
configurations, it might be necessary to get such a control also
further away and this seems to lead to interesting questions,
involving degenerate and non-symmetric Finsler distances.

\subsection{Return to equilibrium, ideas of the proof of Theorem
\ref{kfp2.5}}
\label{kfpre}

 Keeping the same assumptions, let $\Pi_0,\,\Pi _1 $ be the rank 1
spectral projections corresponding to the eigenvalues $\mu _0:=0,\,
\mu _1:=\mu $ of $-\Delta _A^{(0)}$ and put $\Pi =\Pi _0+\Pi _1$. Then
$e_{-1,},e_1$ is a basis for ${\cal R}(\Pi )$ and the restriction of
$P$ to this range, has the matrix \ekv{4.1} {
\left(\begin{array}{ccc}\lambda _{-1}^*\\ \lambda _1^*
\end{array}\right) \left(\begin{array}{ccc} \lambda _{-1} &\lambda
_1\end{array}\right) = \left(\begin{array}{ccc} \lambda _{-1}^*\lambda
_{-1} & \lambda _{-1}^*\lambda _{1}\\ \lambda _{1}^*\lambda _{-1}&
\lambda _{1}^*\lambda _{1}
\end{array}\right) } with the eigenvalues $0$ and $\mu =\lambda
_{-1}^*\lambda _{-1}+ \lambda _{1}^*\lambda _{1}$. A corresponding
basis of eigenvectors is given by
\begin{eqnarray}\label{4.2} v_0&=&\frac{1}{\sqrt{\mu _1}}(\lambda
_1e_{-1}-\lambda _{-1}e_{-1}) \\ v_1&=&\frac{1}{\sqrt{\mu _1}}(\lambda
_{-1}^*e_{-1}+\lambda _{1}^* e_{-1}).\nonumber
\end{eqnarray} The corresponding dual basis of eigenfunctions of $P^*$
is given by \eekv{4.3} { v_0^*&=&\frac{1}{\sqrt{\mu }}(\lambda
_1^*e_{-1}^*- \lambda _{-1}^*e_{-1}^*) }
{v_1^*&=&\frac{1}{\sqrt{\mu}}(\lambda _{-1}e_{-1}^*+\lambda _{1}
e_{1}^*),} where $e_{-1}^*,e_1^* \in {\cal R}(\Pi ^*)$ is the basis
that is dual to $e_{-1},e_1$. It follows that $v_j,v_j^*={\cal O}(1)$
in $L^2$, when $h\to 0$.
\par From this discussion we conclude that $\Pi _j=(\cdot |v_j^*)v_j$,
\emph{are uniformly bounded when $h\to 0$}.  A non-trivial fact, based
on the analysis described in Subsections \ref{kfpgen}, \ref{kfpav}, is
that after replacing the standard norm and scalar product on $L^2$ by
certain uniformly equivalent ones, we have \ekv{4.4} {\Re (Pu|u)\ge
\frac{h}{C}\Vert u\Vert^2,\quad \forall u\in {\cal
R}(1-\widetilde{\Pi} ),} where $\widetilde{\Pi }$ is the spectral
projection corresponding to the spectrum of $P$ in $D(0,Bh)$ for some
$B\gg 1$.

This can be applied to the {study of $u(t):=e^{-tP/h}u(0)$,} where the
initial state $u(0)\in L^2$ is arbitrary: {Write \ekv{4.5} { u(0)=\Pi
_0u(0)+\Pi _1u(0)+(1-\Pi )u(0)=:u^0+u^1+u^\perp .  }Then
\begin{eqnarray}\label{11} \Vert u^0 \Vert, \Vert u^1 \Vert, \Vert
u^\perp\Vert&\le& {\cal O}(1) \Vert u(0)\Vert \\ \Vert
e^{-tP/h}u^\perp\Vert&\le& Ce^{-t/C}\Vert u(0)\Vert \label{12}\\
e^{-tP/h}u_j&=&e^{-t\mu _j/h}u_j,\ j=0,1. \label{13}
\end{eqnarray}} Here (\ref{12}) follows if we write $u^\perp=
(1-\widetilde{\Pi })u+(\Pi -\widetilde{\Pi })u$, apply (\ref{4.4}) to
the evolution of the first term, and use that the last term is the
(bounded) spectral projection of $u$ to a finite dimensional spectral
subspace of $P$, for which the corresponding eigenvalues all have real
part $\ge h/C$.\hfill{$\Box$}

\section{Spectral asymptotics in 2 dimensions}\label{2d}
\setcounter{equation}{0}

\subsection{Introduction}\label{2dint}

This section is mainly based on recent joint works with S.~V{\~ u}
Ng{\d o}c and M.~Hitrik \cite{HiSj} \cite{HiSjVu}, but we shall start
by recalling some earlier results that we obtained with A.~Melin
\cite{MeSj2} where we discovered that in the two dimensional case one
often can have Bohr-Sommerfeld conditions to determine all the
individual eigenvalues in some region of the spectral plane, provided
that we have analyticity. This was first a surprise for us since in
the self-adjoint case such results are known only in 1 dimension and
in very special cases for higher dimensions.

Subsequently, with M. Hitrik we have studied small perturbations of
self-adjoint operators. First we studied the case when the classical
flow of the unperturbed operator is periodic, then also with S.~V{\~
u} Ng{\d o}c we looked at the more general case when it is completely
integrable, or just when the energy surface contains some invariant
diophantine Lagrangian tori.

\subsection{Bohr-Sommerfeld rules in two dimensions}\label{2dbs}

For (pseudo-)differential operators in dimension 1, we often have a
Bohr-Sommerfeld rule to determine the asymptotic behaviour of the
eigenvalues. Consider for instance the semi-classical Schr\"odinger
operator
$$
P=-h^2\frac{d^2}{dx^2}+V(x),\mbox{ with symbol } p(x,\xi )=\xi
^2+V(x),
$$
where we assume that $V\in C^\infty ({\bf R};{\bf R})$ and $V(x)\to
+\infty $, $|x|\to \infty $. Let $E_0\in {\bf R}$ be a non-critical
value of $V$ such that (for simplicity) $\{ x\in {\bf R}; V(x)\le
E_0\}$ is an interval. Then in some small fixed neighborhood of $E_0$
and for $h>0$ small enough, the eigenvalues of $P$ are of the form
$E=E_k$, $k\in {\bf Z}$, where
$$
\frac{I(E)}{2\pi h} = k-\theta (E; h), \quad I(E)=\int_{p^{-1}(E)}\xi
\cdot dx, \ \theta (E;h)\sim \theta _0(E)+\theta _1(E)h+...$$ In the
non-self-adjoint case we get the same results, provided that $\Im V$
is small and $V$ {\it is analytic}. The eigenvalues will then be on a
curve close to the real axis.

\par
For self-adjoint operators in
dimension $\ge 2$ it is generally admitted that Bohr-Sommerfeld rules
do not give all eigenvalues in any fixed domain except in certain
(completely integrable) cases. Using the KAM theorem one can sometimes
describe some fraction of the eigenvalues.

\par With A.~Melin \cite{MeSj2}: we considered an
$h$-pseudodifferential operator with leading symbol $p(x, \xi )$ that
is bounded and holomorphic in a tubular neighborhood of ${\bf R}^4$ in
${\bf C}^4 = {\bf C}^2_x\times{\bf C}^2_\xi $.  Assume that
\ekv{bs.1}{{\bf R}^4 \cap p^{-1}(0)\ne \emptyset \mbox{ is
connected.}}  \ekv{bs.2}{\mbox{On } {\bf R}^4\mbox{ we have }|p(x, \xi
)| \ge 1/C,\mbox{ for }|(x, \xi )| \ge C,} for some $C > 0$,
\ekv{bs.3}{d\Re p(x, \xi ), d \Im p(x, \xi ) \mbox{ are linearly
independent for all } (x, \xi ) \in p^{-1}(0) \cap {\bf R}^4.}  (Here
the boundedness assumption near $\infty $ and (\ref{bs.2}) can be
replaced by a suitable ellipticity assumption.)  It follows that
$p^{-1}(0)\cap {\bf R}^4$ is a compact (2-dimensional) surface.

\par Also assume that \ekv{bs.4}{|\{ \Re p, \Im p\} |\mbox{ is
sufficiently small on } p^{-1}(0) \cap {\bf R}^4.}  Here
``sufficiently small'' refers to some positive bound that can be
defined whenever the the other conditions are satisfied uniformly.

When the Poisson bracket vanishes on $p^{-1}(0)$, this set becomes a
Lagrangian torus, and more generally it is a torus. The following is a
complex version of the KAM theorem without small divisors (cf
T.W.~Cherry \cite{Ch}(1928), J.~Moser \cite{Mo}(1958)),
\begin{theo}\label{bs1} (\cite{MeSj2}) There exists a smooth
2-dimensional torus $\Gamma \subset p^{-1}(0)\cap {\bf C}^4$, close to
$p^{-1}(0)\cap {\bf R}^4$ such that ${{\sigma }_\vert}_{\Gamma } = 0$
and $I_j(\Gamma )\in {\bf R},$ $j = 1, 2$. Here $ I_j(\Gamma ):=
\int_{\gamma _j} \xi \cdot dx$ are the actions along the two
fundamental cycles $\gamma _1,\gamma _2\subset \Gamma $, and $\sigma =
\sum_1^2 d\xi _j \wedge dx_j$ is the complex symplectic (2,0)-form.
\end{theo}

Replacing $p$ by $p-z$ for $z$ in a neighborhood of $0\in {\bf C}$, we
get tori $\Gamma (z)$ depending smoothly on $z$ and a corresponding
smooth action function $I(z)=(I_1(\Gamma (z)),I_2(\Gamma (z)))$, which
are important in the Bohr-Sommerfeld rule for the eigen-values near
$0$ in the semi-classical limit $h\to 0$:
\begin{theo}\label{bs2} (\cite{MeSj2}) Under the above assumptions,
there exists $\theta_0 \in ( \frac{1}{2} {\bf Z})^2$ and $\theta (z;
h) \sim \theta_ 0 + \theta_ 1(z)h + \theta_ 2(z)h^2 + ..$ in $C^\infty
( \mathrm{neigh\,} (0,{\bf C}))$, such that for $z$ in an
$h$-independent neighborhood of 0 and for $h > 0$ sufficiently small,
we have that $z$ is an eigenvalue of $P=p(x,hD_x)$ iff
$$\frac{I(z)}{2\pi h}
= k-\theta (z; h),\mbox{ for some }k ∈\in {\bf Z}^2.\quad (BS)$$
\end{theo}

Recently, a similar result was obtained by S.~Graffi, C.~Villegas Bas
\cite{GrVi}.

\par An application of this result is that we get all resonances
(scattering poles) in a {fixed} neighborhood of $0\in {\bf C}$ for
$-h^2\Delta +V(x)$ if V is an analytic real potential on ${\bf R}^2$
with a nondegenerate saddle point at $x=0$, satisfying $V(0)=0$ and
having $\{ (x,\xi )=(0,0)\}$ as its classically trapped set in the
energy surface $\{ p(x,\xi )=0\}$.

\medskip

\subsection{Diophantine case}\label{2ddi}

In this and the next subsection we describe a result from
\cite{HiSjVu} and the main result of \cite{HiSj} about individual
eigenvalues for small perturbations of a self-adjoint operator with a
completely integrable leading symbol. We start with the case when only
Diophantine tori play a role.

Let $P_\epsilon (x,hD;h)$ on ${\bf R}^2$ have the leading symbol
$p_\epsilon (x,\xi )=p(x,\xi )+i\epsilon q(x,\xi )$ where $p$, $q$ are
real and extend to bounded holomorphic functions on a tubular
neighborhood of ${\bf R}^4$. 
Assume that $p$ fulfills the ellipticity condition (\ref{bs.2})
near infinity and that \ekv{bs.5}{P_{\epsilon =0}=P(x,hD)} is
self-adjoint. (The conditions near infinity can be
modified and we can also replace ${\bf R}^2_x$ by a compact
2-dimensional analytic manifold.)

\par Also, assume that $P_\epsilon (x,\xi ;h)$ depends smoothly on
$0\le \epsilon \le \epsilon _0$ with values in the space of bounded
holomorphic functions in a tubular neighborhood of ${\bf R}^4$, and
$P_\epsilon \sim p_\epsilon +hp_{1,\epsilon} +h^2p_{2,\epsilon }+...$,
when $h\to 0$.

Assume \ekv{bs.6}{p^{-1}(0) \mbox{ is connected and }dp\ne 0 \mbox{ on
that set.}}

Assume complete integrability for $p$: There exists an analytic real
valued function $f$ on $T^*{\bf R}^2$ such that $H_pf = 0$, with the
differentials $df$ and $dp$ being linearly independent almost
everywhere on $p^{-1}(0)$. ($H_p=p'_\xi \cdot \frac{\partial }{\partial x}-p'_x \cdot
\frac{\partial }{\partial \xi }$ is the Hamilton field.)

\par Then we have a disjoint union decomposition
\ekv{bs.7}{p^{-1}(0)\cap T^*{\bf R}^2 = \bigcup_{\Lambda \in J}\Lambda
,} where $\Lambda $ are compact connected sets, invariant under the
$H_p$ flow.  We assume (for simplicity) that $J$ has a natural
structure of a graph whose edges correspond to families of regular
leaves; Lagrangian tori (by the Arnold-Mineur-Liouville theorem).  The
union of edges $J\setminus S$ possesses a natural real analytic
structure.

Each torus $\Lambda \in J\setminus S$ carries real analytic
coordinates $x_1, x_2$ identifying $\Lambda $ with ${\bf T}^2 = {\bf
R}^2/2\pi {\bf Z}^2$, so that along $\Lambda $, we have \ekv{bs.8}{H_p
= a_1\frac{\partial }{\partial x_1} + a_2\frac{\partial }{\partial
x_2},} where $a_1, a_2 \in {\bf R}$. The rotation number is defined as
the ratio $\omega (\Lambda ) = [a_1 : a_2] \in {\bf R}{\bf P}^1$, and
it depends analytically on $\Lambda \in J\setminus S$. We assume that
$\omega (\Lambda )$ is not identically constant on any open edge.

\par We say that $\Lambda \in J\setminus S$ is respectively rational,
irrational, diophantine if $a_1/a_2$ has the corresponding
property. Diophantine means that there exist $\alpha >0$, $d>0$ such
that \ekv{bs.15} { |(a_1,a_2)\cdot k|\ge \frac{\alpha }{|k|^{2+d}},\
0\ne k\in {\bf Z}^2, }

\par We introduce \ekv{bs.9}{\langle q\rangle_T =
\frac{1}{T}\int_{-T/2}^{T/2} q \circ \exp (tH_p) dt,\ T > 0, } and
consider the compact intervals $Q_\infty (\Lambda ) \subset {\bf R}$,
$\Lambda \in J$, defined by, \ekv{bs.10}{Q_\infty (\Lambda ) = [\lim_
{T\to →\infty} \inf_ {\Lambda } \langle q\rangle _T , \lim_ {T\to
\infty} \sup_ {\Lambda} \langle q\rangle _T ] .}

\par { A first localization of the spectrum $\sigma (P_\epsilon
(x,hD_x;h))$ (\cite{HiSjVu}) is given by \ekv{bs.14} { \Im (\sigma
(P_\epsilon ) \cap \{z ; |\Re z|\le \delta \} ) \subset \epsilon [\inf
\bigcup_ {\Lambda \in J} Q_\infty (\Lambda )-o(1), \sup \bigcup_
{\Lambda \in J} Q_\infty (\Lambda ) + o(1) ] , } when $\delta
,\epsilon, h \to 0$.}

\par For each torus $\Lambda \in J\setminus S$, we let {$\langle
q\rangle (\Lambda )$} be the average of ${{q}_\vert}_{\Lambda }$ with
respect to the natural smooth measure on $\Lambda $, and assume that
the analytic function $J\setminus S ∋\ni \Lambda \mapsto \langle
q\rangle (\Lambda )$ is not identically constant on any open edge.

By
combining (\ref{bs.8}) with the Fourier series representation of $q$,
we see that
when $\Lambda$ is irrational then $Q_\infty (\Lambda ) = \{\langle
q\rangle (\Lambda )\}$, while 
in the rational case,    
\ekv{bs.11}{Q_\infty (\Lambda ) \subset \langle q\rangle
(\Lambda ) + {\cal O}( \frac{1} {(|n|+|m|)^\infty}) [ -1, 1 ],}
when $\omega (\Lambda ) = \frac{m}{n}$ and $m \in {\bf Z}$, $n \in {\bf N}$ are relatively prime.

\par{Let $F_0\in \cup_{\Lambda \in J}Q_\infty (\Lambda )$ and assume
that there exists a Diophantine torus $\Lambda _d$ (or finitely many),
such that \ekv{bs.15.5}{ \langle q\rangle (\Lambda _d)=F_0, \quad
d_\Lambda \langle q\rangle (\Lambda _d)\ne 0.}} With M.~Hitrik and
S.~V{\~u} Ng{\d o}c we obtained:

\begin{theo}\label{bs3}(\cite{HiSjVu}) Assume also that $F_0$ does not
belong to $Q_\infty (\Lambda )$ for any other $\Lambda \in J$. Let
$0<\delta <K<\infty $. Then $\exists C>0$ such that for $h>0$ small
enough, and $k^K\le \epsilon \le h^\delta $, the eigenvalues of
$P_\epsilon $ in the rectangle $ |\Re z|<h^\delta /C,\ |\Im z-\epsilon
\Re F_0|<\epsilon h^\delta /C $ are given by
$$P^{(\infty )}(h(k-\frac{k_0}{4})-\frac{S}{2\pi },\epsilon ;h)+{\cal
  O}(h^\infty ),\ k\in {\bf Z}^2,$$ Here $P^{(\infty )}(\xi ,\epsilon
;h)$ is smooth, real-valued for $\epsilon =0$ and when $h\to 0$ we
have \ekv{bs.16}{P^{(\infty )}(\xi ,\epsilon ;h)\sim \sum_{\ell
=0}^\infty h^\ell p_\ell^{(\infty )}(\xi ,\epsilon ),\ p_0^{(\infty
)}=p_\infty (\xi )+i\epsilon \langle q\rangle (\xi )+{\cal O}(\epsilon
^2),} corresponding to action angle coordinates.
\end{theo}

In \cite{HiSjVu} we also considered applications to small
non-self-adjoint perturbations of the Laplacian on a surface of
revolution. Thanks to (\ref{bs.11}) the total measure of the union of
all $Q_\infty (\Lambda )$ over the rational tori is finite and
sometimes small, and we could then show that there are plenty of
values $F_0$, fulfilling the assumptions in the theorem.

\subsection{The case with rational tori}\label{2drt}

{ Let $F_0$ be as in (\ref{bs.15.5}) but now also allow for the
possibility that there is a rational torus (or finitely many) $\Lambda
_r$, such that \ekv{bs.17} { F_0\in Q_\infty (\Lambda _r),\quad
F_0\ne\langle q\rangle (\Lambda _r), } \ekv{bs.17.5} {d_\Lambda
(\langle q\rangle)(\Lambda _r)\ne 0,\ d_\Lambda (\omega )(\Lambda
_r)\ne 0.  } { Assume also that \ekv{bs.18} { F_0\not\in Q_\infty
(\Lambda ),\hbox{ for all }\Lambda \in J\setminus \{ \Lambda _d,
\Lambda _r \} .  }}}
With M.~Hitrik we showed the
following result:

\begin{theo}\label{bs4}(\cite{HiSj}) Let $\delta >0$ be small and
assume that $h\ll \epsilon \le h^{\frac{2}{3}+\delta }$, or that the
subprincipal symbol of $P$ vanishes and that $h^2 \ll \epsilon \le
h^{\frac{2}{3}+\delta }$ . Then the spectrum of $P_\epsilon $ in the
rectangle
$$
[-\frac{\epsilon }{C},\frac{\epsilon }{C}]+i\epsilon
[F_0-\frac{\epsilon ^\delta }{C},F_0+\frac{\epsilon ^\delta }{C}]
$$
is the union of two sets: $E_d \cup E_r$, where the elements of $E_d$
form a distorted lattice, given by the Bohr-Sommerfeld rule
(\ref{bs.16}), with horizontal spacing $\asymp h$ and vertical spacing
$\asymp \epsilon h$. The number of elements $\# (E_r)$ of $E_r$ is
${\cal O}({\epsilon ^{3/2}}/{h^2})$.
\end{theo} NB that $\# (E_d) \asymp \epsilon ^{1+\delta }/h^2$.

This result can be applied to the damped wave equation on surfaces of
revolution.

\subsection{Outline of the proofs of Theorem \ref{bs3} and \ref{bs4}}
\label{2dpr}

The principal symbol of $P_\epsilon $ is $p_\epsilon =p+i\epsilon
q+{\cal O}(\epsilon ^2)$. Put
$$
\langle q\rangle_T=\frac{1}{T}\int_{-T/2}^{T/2}q\circ \exp (tH_p)dt.
$$
As in Section \ref{kfp} we will use an averaging of the imaginary part
of the symbol.  Let $J(t)$ be the piecewise affine function with
support in $[-\frac{1}{2},\frac{1}{2}]$, solving
$$
J'(t)=\delta (t)-1_{[-\frac{1}{2},\frac{1}{2}]}(t),
$$
and introduce the weight
$$G_T(t)=\int J(-\frac{t}{T})q\circ \exp (tH_p)dt.$$
Then $H_pG_T=q-\langle q\rangle_T$, implying \ekv{2dpr.1} {p_\epsilon
\circ \exp (i\epsilon H_{G_T})=p+i\epsilon \langle q\rangle_T+{\cal
O}_T(\epsilon ^2).}

The left hand side of (\ref{2dpr.1}) is the principal symbol of the
isospectral operator $e^{-\frac{\epsilon }{h}G_T(x,hD_x)}\circ
P_\epsilon \circ e^{\frac{\epsilon }{h}G_T(x,hD_x)}$ and under the
assumptions of Theorem \ref{bs3} resp. \ref{bs4} its imaginary part
will not take the value $i\epsilon F_0$ on $p^{-1}(0)$ away from
$\Lambda _d$ resp.~$\Lambda _d\cup \Lambda _r$. This means that we
have localized the spectral problem to a neighborhood of $\Lambda _d$
resp.~$\Lambda _d\cup\Lambda _r$.

Near $\Lambda _d$ we choose action-angle coordinates so that $\Lambda
_d$ becomes the zero section in the cotangent space of the 2-torus,
and \ekv{2dpr.2} {p_\epsilon (x,\xi )=p(\xi )+i\epsilon q(x,\xi
)+{\cal O}(\epsilon ^2).}  We follow the quantized Birkhoff normal
form procedure in the spirit of V.F.~Lazutkin and Y.~Colin de
Verdi\`ere \cite{La, Co}: solve first \ekv{2dpr.3}{H_pG=q(x,\xi
)-\langle q(\cdot ,\xi )\rangle ,} where the bracket indicates that we
take the average over the torus with respect to $x$. Composing with
the corresponding complex canonical transformation, we get the new
conjugated symbol
$$p(\xi )+i\epsilon \langle q(\cdot ,\xi )\rangle+{\cal O}(\epsilon
^2+\xi ^\infty ).$$ Here the Diophanticity condition is of course
important.

Iterating the procedure we get for every $N$,
$$
p_\epsilon \circ \exp (H_{G^{(N)}})=\underbrace{p(\xi )+i\epsilon
(\langle q\rangle (\xi )+{\cal O}(\epsilon ,\xi
))}_{\mathrm{independent\ of\ }x}+{\cal O}((\xi ,\epsilon )^{N+1})
$$
This procedure can be continued on the operator level, and up to a
small error we see that $P_\epsilon $ is microlocally equivalent to an
operator $P_\epsilon (hD_\xi ,\epsilon ;h)$. At least formally,
Theorem \ref{bs3} then follows by considering Fourier series
expansions, but in order to get a full proof we also have take into
account that we have constructed complex canonical transformations
that are quantized by Fourier integral operators with complex phase
and study the action of these operators on suitable exponentially
weighted spaces.

Near $\Lambda _r$ we can still use action-angle coordinates as in
(\ref{2dpr.2}) but the homological equation (\ref{2dpr.3}) is no
longer solvable. Instead, we use {secular perturbation theory} (cf the
book \cite{LiLi}), which amounts to making a {partial Birkhoff
reduction.}

\par After a linear change of $x$-variables, we may assume that $p(\xi
)=\xi _2+{\cal O}(\xi ^2)$ and in order to fix the ideas $=\xi _2+\xi
_1^2$. Then we can make the averaging procedure only in the
$x_2$-direction and reduce $p_\epsilon $ in (\ref{2dpr.2}) to
$$
\widetilde{p}_\epsilon (x,\xi )=\underbrace{\xi _2+\xi _1^2+{\cal
O}(\epsilon )}_{\mathrm{independent\ of\ }x_2,\atop \approx \xi _2+\xi
_1^2+i\epsilon \langle q\rangle_2(x_1,\xi )}+{\cal O} ((\epsilon ,\xi
)^\infty ),
$$
where $\langle q\rangle_2(x_1,\xi )$ denotes the average with respect
to $x_2$.

\par
Carrying out the reduction on the operator level, we obtain up to
small errors an operator $\widetilde{P}_\epsilon
(x_1,hD_{x_1},hD_{x_2};h)$ and after passing to Fourier series in
$x_2$, a family of non-self-adjoint operators on $S_{x_1}^1$:
$\widetilde{P}_\epsilon (x_1,hD_{x_1},hk;h)$, $k\in {\bf Z}$.

The non-self-adjointness and the corresponding possible wild growth of
the resolvent makes it hard to go all the way to study individual
eigenvalues. However, it can be shown that in the region $|\xi _1|\gg
\epsilon ^{1/2}$ (inside the energy surface $p=0$) we can go further
and (as near $\Lambda _d$) get a sufficiently good elimination of the
$x$-dependence. This leads to the conclusion that the contributions
from a vicinity of $\Lambda _r$ to the spectrum of $P_\epsilon $ in
the rectangle
$$
|\Re z|\le \frac{\epsilon }{C},\ |\Im z-\epsilon F_0|\le
\frac{\epsilon ^{1+\delta }}{C},
$$ 
come from a neighborhood of $\Lambda _r$ of phase space volume ${\cal
O}(\epsilon ^{3/2})$.

This explains heuristically why the rational torus will contribute
with ${\cal O} (\epsilon ^{3/2}/h^2)$ eigenvalues in the rectangle.

The actual proof is more complicated. We use a Grushin problem
reduction in order to reduce the study near $\Lambda _r$ to that of a
square matrix of size ${\cal O}(\epsilon ^{3/2}/h^2)$. However, even
if we avoid the eigenvalues of such a matrix, the inverse can only be
bounded by \ekv{2dpr.4}{ \exp {\cal O}(\epsilon ^{3/2}/h^2).} What
saves us is that away from $\Lambda _r\cup \Lambda _d$, we can
conjugate the operator with exponential weights and show that the
resolvent has an ``off-diagonal decay'' like $\exp (-1/(Ch))$. This
implies that we can confine the growth in (\ref{2dpr.4}) to a small
neighborhood of $\Lambda _r$, if
$$
\frac{1}{Ch}\gg \frac{\epsilon ^{\frac{3}{2}}}{h^2},
$$ 
leading to the assumption $\epsilon \ll h^{2/3}$ in Theorem \ref{bs4}.

\section{Weyl asymptotics for non-self-adjoint operators}\label{r}
\setcounter{equation}{0}
\subsection{Introduction}

\par For self-adjoint differential (pseudo)differential operators we
have (under suitable assumptions) the Weyl law for the asymptotic
distribution of eigenvalues, established in higher dimensions by H.
Weyl \cite{We} in 1912 in the case of second order elliptic boundary
value problems.

In the semiclassical setting such results were obtained by
J. Chazarain, B.Helffer--D.Robert, V.Ivrii and many others (see
\cite{DiSj} and further references there).  Under suitable additional
assumptions it states that if $P=P^w(x,hD_x;h)$ is a self-adjoint
$h$-pseudodifferential operator with leading (real) symbol $p(x,\xi
)$, then if $\Omega \subset{\bf C}$ is a domain intersecting ${\bf R}$
along a bounded interval, the number of eigenvalues of $P$ in $\Omega
$ (intersected with ${\bf R}$) satisfies
$$
\#(\sigma (P)\cap\Omega )={1\over (2\pi h)^n}({\rm vol\,}(p^{-1}
(\Omega ))+o(1)),\ h\to 0.
$$

A simple example is provided by the semiclassical harmonic oscillator
${1\over 2}((hD_x)^2+x^2)$ on the real line which has the eigenvalues
$(k+{1\over 2})h$, $k=0,1,...$.

\par In the non-self-adjoint case Weyl-asymptotics is known to hold in
some cases close to the self-adjoint case or for normal operators.
 \par We do not always have Weyl-asymptotics: Following Davies and
Boulton (see \cite{Da2}), we can consider the non-self-adjoint
harmonic operator: $P={1\over 2}((hD_x)^2+ix^2)$ whose eigenvalues are
given by $e^{i\pi /4}(k+{1\over 2})h$, $k\in{\bf N}$ (\cite{Sj}):

The set of values of $p={1\over 2}(\xi ^2+ix^2)$ is the closed first
quadrant and if we choose the open bounded set $\Omega $ to intersect
the 1st quadrant but not the line $\Im z=\Re z $, we get ${\rm
vol\,}(p^{-1} (\Omega ))>0$, while there are no eigenvalues in $\Omega
$.

\par More generally, $h$-differential operators with analytic
coefficients often have their spectrum determined by complex-geometric
quantities, and are likely not to obey the Weyl law.  In particular,
in the one dimensional case it often happens that the eigenvalues are
concentrated to certain curves with branch points.

\par As we have seen in Theorem \ref{int3} we are often confronted
with the pseudospectral phenomenon: On the image of $p$ the resolvent
may be very large even far from the spectrum. This causes the
eigenvalues to be very sensitive to small perturbations of the
operator (by Theorem \ref{int2}).

\par In her thesis M. Hager (see \cite{Ha2}) considered a class of
perturbed $h$-pseudodifferential operators on the real line of the
form $P_\delta =P(x,hD;h)+\delta q_\omega (x)$, where $P$ is analytic
and $q_\omega $ is a random linear combination of the $C/h$ first
eigen-functions of an auxiliary operator. She showed that with
probability very close to 1 when $h\to 0$, $P_\delta $ obeys Weyl
asymptotics.

\par Here, we shall discuss a generalization to the multidimensional
case obtained with Hager \cite{HaSj}.  The results will be much more
general in many ways, but the class of perturbations will be slightly
different.
 
\subsection{The result}

\par\noindent a) {\it The unperturbed operator.}  Let $m(\rho )\ge 1$,
$\rho =(x,\xi )$ be an order function on ${\bf R}^{2n}$, so
that $$0<m(\rho )\le C_0\langle \rho -\mu \rangle ^{N_0}m(\mu ),$$ we
may assume that $m\in S(m)=$ $\{ u\in C^\infty ({\bf
R}^{2n});\, \partial _\rho ^\alpha u={\cal O}_\alpha (m),\,\forall
\alpha \in{\bf N}^{2n}\}$. Assume $m\ge 1$ and let
$$
P(\rho ;h)\sim p(\rho )+hp_1(\rho )+... \hbox{ in } S(m).
$$ 
Assume $\exists z_0\in{\bf C}$, $C_0>0$ such that
$$
\vert p(\rho )-z_0\vert \ge m(\rho )/C_0 \hbox{ (ellipticity)}.
$$
Let
$$
\Sigma =\overline{p({\bf R}^{2n})}=p({\bf R}^{2n})\cup\Sigma _\infty ,
$$
$$
\Sigma _\infty =\{ \lim_{j\to\infty }p(\rho _j);\ {\bf R}^{2n}\ni\rho
_j\to \infty \}
$$
We write: $P=P^w(x,hD;h)$.

\par Let $\Omega \subset\subset {\bf C}\setminus \Sigma _\infty $ be
open and simply connected containing $z_0$. Then using the
pseudodifferential calculus, it is easy to show: \smallskip
\par\noindent 1) $\sigma (P)\cap\Omega $ is discrete when $h>0$ is
small enough.  \smallskip
\par\noindent 2) $\forall \epsilon >0$, $\exists h(\epsilon )>0$, such
that $\sigma (P)\cap \Omega \subset \Sigma +D(0,\epsilon )$, $0<h\le
h(\epsilon ).$

\par
\medskip
\par\noindent {b) \it The random pertubation. } Let
$0<\widetilde{m},\widehat{m}\le 1$ be square integrable order
functions, one of which is integrable. Let $\widetilde{S}\in
S(\widetilde{m})$, $\widehat{S}\in S(\widehat{m})$ be elliptic
symbols. The corresponding operators are Hilbert-Schmidt with $\Vert
\widetilde{S}\Vert _{\rm HS},\Vert \widehat{S}\Vert _{{\rm HS}}={\cal
O}(h^{-n/2})$. Let
$$
Q_\omega =\widehat{S}\circ \sum_{j,k}\alpha _{j,k}(\omega
)\widehat{e}_j\widetilde{e}_k^*\circ \widetilde{S},
$$
where $\alpha _{j,k}$, $j,k\in{\bf N}$, are independent complex
Gaussian random variables with expectation value 0 and variance $1$,
and $(\widehat{e}_j)_1^\infty $ and $(\widetilde{e}_j)_1^\infty $ are
orthonormal bases in $L^2({\bf R}^n)$,
$\widehat{e}_j\widetilde{e}_k^*u=(u\vert
\widetilde{e}_k)\widehat{e}_j$.

\par Let $M=C_1h^{-n}$ with $C_1\gg 1$. Then with probability $\ge
1-Ce^{-h^{-2n}/C}$ the Hilbert-Schmidt and trace class norms of $Q$ fulfil: 
\ekv{r1.1} {\Vert Q\Vert _{\rm HS}\le M,\ \Vert
Q\Vert _{\rm tr}\le M^{3\over 2}}

Let $\Gamma \subset\subset \Omega $ be open with smooth
boundary. \medskip
\begin{theo}\label{r1} (\cite{HaSj}).  Assume 
\ekv{r1.2}{ p(\rho
)\in\partial \Gamma \Rightarrow dp(\rho ),\,d\overline{p}(\rho )
\hbox{ are linearly independent.}}

Let $\epsilon ,\delta $ depend on $h$ with $0<\epsilon \ll 1$,
$$
e^{-{\epsilon\over Ch }}\le \delta \ll h^{3n+\frac{1}{2}},\ C\gg 1
$$ 
(implying that $\epsilon \ge {\rm Const.}h\ln{1\over h}$). Then with
probability $\ge 1-{C\over \sqrt{\epsilon }}e^{-{\epsilon\over 2(2\pi
h)^n }}$, we have \ekv{r1.3} { \vert \# (\sigma (P_\delta )\cap\Gamma
)-{1\over (2\pi h)^n}{\rm vol\,}(p^{-1} (\Gamma ))\vert \le
C{\sqrt{\epsilon }\over h^n}.  .}\end{theo}

\par There is a similar result giving (\ref{r1.3}) simultaneously for
all $\Gamma $ in a suitable family.

The assumption (\ref{r1.2}) implies that $\partial \Gamma \cap\partial
\Sigma =\emptyset$, so if we want to count the eigenvalues near
$\partial \Sigma$, we need to weaken that assumption.
For $z\in{\rm neigh\,}(\partial \Gamma )$, we put 
\ekv{r1.4} {
V_z(t)={\rm vol\,}\{ \rho \in{\bf R}^{2n};\, \vert p(\rho )-z\vert
^2\le t\} .  }
Introduce the assumption 
\ekv{r1.5} { \exists \kappa\in ]0,1],\hbox{
such that }V_z(t)={\cal O}(t^{\kappa}),\hbox{ uniformly for
}z\in{\rm neigh\,}(\partial \Gamma ),\ 0\le t\ll 1.}
\par\noindent {\it Example. } (\ref{r1.2}) $\Rightarrow$ (\ref{r1.5})
with $\kappa =1$.  \smallskip
\par\noindent {\it Example. } The best that can happen when $\partial
\Gamma \cap\partial \Sigma \ne \emptyset$ is that \ekv{r1.6} { p(\rho
)\in\partial \Gamma \Rightarrow \{ p,\overline{p}\}(\rho )\ne 0\hbox{
or }\{ p,\{ p,\overline{p}\} \}(\rho )\ne 0.  } It is easy to see that
(\ref{r1.6}) implies (\ref{r1.5}) with $\kappa =3/4$. (\ref{r1.6})
holds for the non-self-adjoint harmonic oscillator when
$0\not\in\partial \Gamma $.

\begin{theo}\label{r2}(\cite{HaSj}) We assume (\ref{r1.5}). Let $\epsilon ,\delta $
depend on $h$ with $0<\epsilon \ll 1$,
$$
e^{-{\epsilon \over Ch^{\kappa }}}\le \delta \ll h^{3n+\frac{1}{2}},\
C\gg 1,
$$
implying that $\epsilon \ge {\rm Const.}h^{\kappa }\ln{1\over
h}$. Then for $0<r\ll 1$ we have with probability $\ge 1-{C\over
r}e^{-{\epsilon \over 2}(2\pi h)^{-n}}$, that \eekv{r1.7} { \vert
\#(\sigma (P_\delta )\cap\Gamma )-{1\over (2\pi h)^n}{\rm
vol\,}(p^{-1} (\Gamma ))\vert \le } {{C\over h^n}\left( {\epsilon
\over r}+C_N\left( r^N+\ln({1\over r}) {\rm vol\,}\left( p^{-1}
(\partial \Gamma +D(0,r))\right)\right)\right) , } for every fixed
$N\in{\bf N}$.\end{theo}

\par If $\kappa >{1\over 2}$, we have ${\rm vol\,}(p^{-1} (\partial
\Gamma +D(0,r)))={\cal O}(r^{2\kappa -1})$ with $2\kappa -1>0$ and in
all cases {we may assume} that $\ln ({1\over r}){\rm vol\,}(p^{-1}
(\partial \Gamma +D(0,r)))={\cal O}(r^{\alpha _0})$, where $\alpha
_0>0$.  Then we choose $N\gg 1$, $r=\epsilon ^{1/(1+\alpha _0)}$ and
the right hand side of (\ref{r1.7}) becomes ${\cal O}(1)\epsilon
^{\alpha _0/(1+\alpha _0)}h^{-n}$.

\par Again we have a similar theorem where the conclusion (\ref{r1.7})
is valid simultaneously for all $\Gamma $ in a suitable family.

\par Recently, the author obtained similar results when $Q_\omega $ is
an operator of multiplication, see \cite{Sj4}.

\subsection{Outline of the proofs}

\par We can construct $\widetilde{P}\in S(m)$ such that $\widetilde{P}(\rho
;h)=P(\rho ;h)$ for $\vert \rho \vert \gg 1$ and $\vert
\widetilde{P}(\rho ;h)-z\vert \ge m(\rho )/C$ for $\rho \in{\bf
R}^{2n}$, $z\in\overline{\Omega }$. The eigenvalues of $P$ in $\Omega
$ coincide with the zeros of the holomorphic function \eekv{r2.1}
{F(z;h)&=&\det P_z,}{P_z&=&(\widetilde{P}(x,hD;h)-z)^{-1}
(P(x,hD;h)-z)} The same remark holds for $P_\delta $ and $F_\delta $
defined as in (\ref{r2.1}) with $P$ replaced by $P_\delta $, provided
that (\ref{r1.1}) holds.

\par For $z\in{\rm neigh\,}(\partial \Gamma )$, put $Q=P_z^*P_z$. Let
$1_\alpha (E)=\max (E,\alpha )$, where $\alpha =Ch$, $C\gg 1$. Using
semiclassical analysis, we can show that under the assumption
(\ref{r1.5}) (cf \cite{MeSj}) \ekv{r2.2} { \ln \det Q\le \ln \det
1_\alpha (Q)= {1\over (2\pi h)^n}(\iint \ln q dxd\xi +{\cal
O}(1)h^{\kappa } \ln{1\over h }), } where $q=\vert p_z\vert ^2$ is the
leading symbol of $Q$.

\par Since
$$
\ln\det Q=\ln \det P_z^*P_z=\ln \vert F(z;h)\vert ^2
$$
we have \ekv{r2.3} { \ln \vert F(z;h)\vert \le {1\over (2\pi
h)^n}(\iint \ln \vert p_z\vert dxd\xi +{\cal O}(1) h^{\kappa }
\ln{1\over h }), } For $\delta >0$ small enough, we
get the same upper bound for $\ln \vert F_\delta (z;h)\vert $
(provided that (\ref{r1.1}) holds).

\par The {main step in the proof} is to get a corresponding lower
bound for each fixed $z$ with a probability close to 1. In the
multidimensional case this boils down to a question about random
determinants. Let $z\in{\rm neigh\,}(\partial \Gamma )$. Let $e_1,e_2,
...$ be the first eigen-functions of $Q=P_z^*P_z$ and let $f_1,f_2,...$ be
the first eigen-functions for $P_zP_z^*$. The two operators have the
same eigenvalues $0\le\lambda _1\le \lambda _2\le ..$.

We can arrange so that
$$
P_ze_j=\sqrt{\lambda _j}f_j,\ P_z^*f_j=\sqrt{\lambda _j}e_j.
$$
Let $N=N(\alpha )=\#\{ j;\lambda _j\le \alpha \}$ ($\alpha =Ch$, $C\gg
1$). Semiclassical analysis gives that
$$
N={\cal O}(h^{\kappa -n}).
$$

\par Consider
$$
{\cal P}^0= \left(\begin{array}{ccc} P_z &R_-\\ R_+
&0 \end{array}\right): L^2\times {\bf C}^N\to L^2\times {\bf C}^N,
$$
$$
R_+:L^2\to {\bf C}^N,\ R_+u(j)=\sqrt{\alpha }(u\vert e_j),
$$
$$
R_-:{\bf C}^N\to L^2,\ R_-=\sqrt{\alpha }\sum_1^N u_-(j)f_j.
$$

${\cal P}^0$ has an inverse
$$
{\cal E}^0=\left(\begin{array}{ccc}E^0 &E_+^0\cr E_-^0
&E_{-+}^0 \end{array}\right) ={\cal O}({1\over \sqrt{\alpha }}),
$$
with $E_+^0,E_-^0,E_{-+}^0$ "explicit", and
$$
\ln \vert \det {\cal P}^0\vert ^2=N\ln \alpha +\det 1_\alpha
(P_z^*P_z).
$$

\par For $P_\delta =P+\delta Q_\omega $, we form
$$
P_z^\delta =(\widetilde{P}-z)^{-1} (P-z+\delta Q_\omega )=P_z+\delta
\widetilde{Q}_\omega ,
$$ and since $\delta \Vert Q_\omega \Vert \le \delta C_0 h^{-n}\ll 1$, $$
{\cal P}^\delta :=\left(\begin{array}{ccc} P_z^\delta &R_-\cr R_+
&0\end{array}\right) $$ is invertible with inverse $$ {\cal E}^\delta
=\left(\begin{array}{ccc} E^\delta &E_+^\delta \cr E_-^\delta
&E_{-+}^\delta \end{array}\right)\approx {\cal E}^0.
$$

Here \ekv{r2.4}{E_{-+}^\delta =E_{-+}^0+\delta
E_-^0\widetilde{Q}_\omega E_+^0+{\rm "small"},} and we can show by
perturbative arguments that
$$
\ln \det {\cal P}^\delta =\ln \det {\cal P}^0+{\cal O}({\delta \over
\sqrt{\alpha }}M^{3/2}),
$$ 
leading to \ekv{r2.5} {\hskip -1truecm \ln\vert \det {\cal P}^\delta
\vert ={1\over (2\pi h)^n}(\iint \ln \vert p_z\vert dxd\xi +{\cal
O}(h^{\kappa }\ln{1\over h})).}

\par On the other hand, computations in \cite{SjZw2} can be used to
get \ekv{r2.6} { \ln \vert \det P_z^\delta \vert =\ln \vert \det {\cal
P}^\delta \vert +\ln \vert \det E_{-+}^\delta \vert } Using
(\ref{r2.4}), we can view $E_{-+}^\delta $ as a random matrix of size
${\cal O}(h^{\kappa -n})$, close in a suitable sense to one with
independent Gaussian random variables as its entries. This can be used
to show:

\par {For every $z\in{\rm neigh\,}(\partial \Gamma )$, we have a nice
lower bound on $\ln \vert \det E_{-+}^\delta \vert $ with probability
close to 1.} (\ref{r2.6}) then gives a corresponding lower bound on
$\ln \vert \det P_z^\delta \vert $.

 To complete the proof of Theorem \ref{r1} we can apply the following
result of M. Hager \cite{Ha1, Ha2} with $\widetilde{h}=h^n$:

\begin{prop} (\cite{Ha1, Ha2}). \it Let $\Gamma $ and $\Omega $ be as
above. Let $\phi \in C(\Omega ;{\bf R})$ be smooth near $\partial
\Gamma $. Let $f=f(z;\widetilde{h})$ be holomorphic in $\Omega $ with
$$
\vert f(z;\widetilde{h})\vert \le e^{\phi (z)/\widetilde{h}},\
z\in{\rm neigh\,}(\partial \Gamma ),\ 0<\widetilde{h}\ll 1.
$$
Assume there exist $\epsilon =\epsilon (\widetilde{h}) \ll 1$,
$z_k=z_k(\widetilde{h})\in\Omega $, $k\in J=J(h)$, such that
\begin{eqnarray*} &&\partial \Gamma \subset \bigcup_{k\in
J}D(z_k,\sqrt{\epsilon }),\ \# J={\cal O}({1\over \sqrt{\epsilon
}}),\\ && \vert f(z_k;\widetilde{h})\vert \ge e^{ (\phi (z_k)-\epsilon
)/\widetilde{h}},\ k\in J.
\end{eqnarray*} Then,
\begin{eqnarray*} \# (f^{-1}(0)\cap \Gamma )={1\over 2\pi
\widetilde{h}}\iint_\Gamma (\Delta \phi )d(\Re z)d(\Im z)+{\cal
O}({\sqrt{\epsilon }\over \widetilde{h}}).
\end{eqnarray*}\end{prop}

\par For the proof of Theorem 1.2 we use an improved version of this
result, see \cite{HaSj}.

\subsection{Comparison with Theorem \ref{bs2}}

From the example with the non-self-adjoint harmonic oscillator in
dimension 1, we have seen that Weyl asymptotics does not always hold
for differential operators in one dimension, when the coefficients are
analytic.  If we add a small random perturbation to the
non-self-adjoint harmonic oscillator, the theorems above and the main
result in \cite{Ha2} show that with probability close to 1 the
eigenvalues will no longer be confined to a half-line but will tend to
fill up the range of the principal symbol $p$ with a density that is
given by $(2\pi h)^{-n}p_*(dv )$, where $dv$ denotes the symplectic
volume element on ${\bf R}^{2n}$ and $p_*(dv)$ is the direct image
under $p$.

From this simple one-dimensional example it is easy to build examples
in higher dimension when Weyl asymptotics does not hold.  In the
2-dimensional case, we can also consider the situation when the
unperturbed operator $P$ satisfies the assumptions of Theorem
\ref{bs2}. It is then natural to compare the distribution law given by
Theorem \ref{bs2} for $P$ and the one given by the theorems \ref{r1},
\ref{r2} for the random perturbations. To leading order in $h$, we get
Weyl asymptotics already for $P$ in the (close to normal) case when
$\{ p, \overline{p}\}$ vanishes identically. In general however, we
get different asymptotic distributions already to leading order
(\cite{Sj3}).

\end{document}